\documentclass[11pt]{amsart}
\usepackage{amstext}
\usepackage{amsthm}
\usepackage{amsopn}
\usepackage{amsfonts}
\usepackage{amsmath}
\usepackage{latexsym}
\usepackage{url}
\usepackage{amscd}
\usepackage{amssymb}
\usepackage{amsmath}
\usepackage{xypic}
\xyoption{all}
\swapnumbers
\newtheorem{thm}[equation]{Theorem}
\newtheorem{prop}[equation]{Proposition}
\newtheorem{lem}[equation]{Lemma}
\newtheorem{cor}[equation]{Corollary}

\theoremstyle{definition}

\newtheorem{remark}[equation]{Remark}
\newtheorem{example}[equation]{Example}

\numberwithin{equation}{section}

\def\limind{\mathop{\oalign{lim\cr
\hidewidth$\longrightarrow$\hidewidth\cr}}}

\newcommand{\Ima}{\operatorname{Im}}

\newcommand{\bbG}{{\mathbb G}}

\newcommand{\bbZ}{{\mathbb Z}}

\newcommand{\Spec}{\operatorname{Spec}}

\newcommand{\GL}{{\operatorname{GL}}}
\newcommand{\SL}{{\operatorname{SL}}}

\newcommand{\GO}{{\operatorname{GO}}}
\newcommand{\Orth}{{\operatorname{O}}}

\newcommand{\Spin}{{\operatorname{Spin}}}

\newcommand{\oG}{\overline{G}}

\newcommand{\rank}{\operatorname{rank}}
\newcommand{\Stab}{\operatorname{Stab}}

\newcommand{\Char}{\operatorname{char}} 

\newcommand{\cd}{\operatorname{cd}}
\newcommand{\ed}{\operatorname{ed}}

\newcommand{\Gal}{\operatorname{Gal}}

\newcommand{\trdeg}{\operatorname{trdeg}}

\begin{document}

\title[Lower bounds]{A lower bound on the
essential dimension of  a connected linear group}
\author{P. Gille}
\address{UMR 8552 du CNRS, DMA, Ecole Normale
Sup\'erieure, F-75005 Paris, France}
\email{Philippe.Gille@ens.fr}

\author{Z. Reichstein}
\address{Department of Mathematics, University of British Columbia,
 Vancouver, BC V6T 1Z2, Canada} \email{reichst\@@math.ubc.ca}
\thanks{Z. Reichstein was partially supported by an NSERC research
grant.}

\subjclass{11E72, 20G10, 14L30}

\keywords{Linear algebraic group, essential dimension,
non-abelian cohomology, group action}

\begin{abstract} Let $G$ be a connected linear algebraic
group defined over an algebraically closed field $k$ and
$H$ be a finite abelian subgroup of $G$ whose order is not 
divisible by $\Char(k)$. We show that the essential dimension 
of $G$ is bounded from below by $\rank(H) - \rank \, C_G(H)^0$,
where $\rank \, C_G(H)^0$ denotes the rank of the maximal torus in
the centralizer $C_G(H)$. This inequality, conjectured by J.-P. Serre,
generalizes previous results of Reichstein -- Youssin 
(where $\Char(k)$ is assumed to be $0$ and $C_G(H)$ 
to be finite) and Chernousov -- Serre (where $H$ is 
assumed to be a $2$-group).
\end{abstract}


\keywords{Central simple algebra, algebraic group,
group action, geometric action, tame action, semisimple
algebra, Galois cohomology}

\maketitle
\tableofcontents

\section{Introduction}

Let $k$ be a base field, $K/k$ be a field extension, $G/k$ be
a linear algebraic group
and $\alpha \in H^1(K, G)$ be a $G$-torsor over $\Spec(K)$.
We will say that $\alpha$ {\em descends} to a subfield $K_0 \subset K$
if $\alpha$ lies in the image of the natural map
$H^1(K_0, G) \to H^1(K, G)$.
The essential dimension $\ed_k(\alpha)$ of $\alpha$,
is defined as the minimal value of $\trdeg_k(K_0)$,
where $\alpha$ descends to $K_0$ and $k \subset K_0$.
(Throughout this paper we will work over a fixed algebraically 
closed field $k$; for this reason we will write $\ed(\alpha)$ 
in place of $ed_k$.)
We also define $\ed(\alpha; l)$ as the minimal value of $\ed(\alpha_L)$,
as $L/K$ ranges over all finite field extensions of degree prime to $l$.
Here $l$ is a prime integer.  The essential dimension 
$\ed(G)$ of the group $G$ (respectively, the essential
dimension $\ed(G; l)$ of $G$ at $l$) is defined 
as the maximal value of $\ed(\alpha)$ (respectively, of
$\ed(\alpha; l)$), as $K/k$  ranges over all field extensions 
and $\alpha$ ranges over $H^1(K, G)$.

Note that the numbers $\ed(\alpha)$, $\ed(\alpha; l)$, $\ed(G)$ and
$\ed(G; l)$ all depend on the base field $k$, which will be
assimed to be fixed (and algebraically closed) throughout this paper.
For details on the notion of essential dimension,
its various interpretations and numerous examples,
see~\cite{reichstein1},~\cite{ry} and~\cite{bf}.
Many of the best known lower bounds on $\ed(G)$
and $\ed(G; l)$, for specific groups $G$, are
deduced from the following inequality. 

\begin{thm} \label{thm.ry} {\rm (}\cite[Theorem 7.8]{ry}{\rm )}
Let $G$ be a connected semisimple linear algebraic group defined
over an algebraically  closed field $k$ of characteristic zero and
$H$ be a finite abelian
subgroup of $G$. Assume that the centralizer $C_G(H)$ is finite.
Then

\smallskip
(a)  $\ed(G) \ge \rank(H)$ and

\smallskip
(b) if $H$ is a $l$-group then $\ed(G; l) \ge \rank(H)$.
\end{thm}

Here by the rank of a finite abelian group $H$ we mean the smallest positive
integer $r$ such that $H$ can be written as a direct product of $r$
cyclic groups.  Equivalently, $\rank(H)$ is the minimal dimension of
a faithful complex linear representation of $H$.

The purpose of this paper is to prove the following more general inequality
conjectured by J.-P. Serre (e-mail, July 25, 2005).

\begin{thm} \label{thm.main}
Let $G$ be a connected reductive linear algebraic group defined
over an algebraically closed base field $k$. Suppose that $H$
is a finite abelian subgroup of $G$
and $\Char(k)$ does not divide $|H|$. Then

\smallskip
(a) $\ed(G) \ge \rank \, H - \rank \, C_G(H)^0$.

\smallskip
(b) Moreover, if $H$ is a $l$-group then
$\ed(G; l) \ge \rank \, H - \rank \, C_G(H)^0$.
\end{thm}

Here $C_G(H)^0$ denotes the connected component of
the centralizer of $H$ in $G$, and by the rank of
this connected group we mean the dimension of its maximal
torus.  In particular, if $\Char(k) = 0$ and
the centralizer $C_G(H)$ is finite
(i.e., $\rank \, C_G(H)^0 = 0$) then
Theorem~\ref{thm.main} reduces
to Theorem~\ref{thm.ry}. Note, however, that even in this
special case the proof we present here is simpler than the one
in~\cite{ry}; in particular, it does not rely on resolution
of singularities.

We also remark that our argument shows a bit more, namely that
the essential dimension of a particular torsor, which we call
a {\em loop torsor}, is $\ge \rank(H) - \rank \, C_G(H)^0$.
Here by a loop torsor we mean the image of a versal $H$-torsor under
the natural map $H^1(\, * \, , H) \to H^1(\, * \,, G)$. (Such torsors
come up in connection with loop algebras; see~\cite{GP}.)

Chernousov and Serre~\cite{cs} used techniques from the theory
of quadratic forms to show that, in the case where $H$ is
a $2$-group, many of the bounds given by Theorem~\ref{thm.ry}(b)
remain valid over any algebraically
closed field base field $k$ of characteristic $\ne 2$.
The ``incompressible" quadratic forms they construct are
closely related to loop torsors; our arguments may
thus be viewed as extending their approach to abelian
subgroups $H$ which are not necessarily $2$-groups.

In order to clarify the exposition we will give two proofs
of Theorem~\ref{thm.main}. The first one, presented in
Section~\ref{sect.char0}, is quite short but it relies
on resolution of singularities and, in particular,
only works in characteristic zero.  The second proof,
presented in Section~\ref{char.free(a)} requires
a bit more work. The advantage of this more elaborate argument
is that it is entirely independent of resolution
of singularities; in particular, it works in prime
characteristic not dividing $|H|$.
Both proofs rely, in a key way, on the existence results
for wonderful (and regular) group compactifications
from~\cite{brion1} and~\cite{brion-kumar} (see
Section~\ref{sect.wonderful}) and on the ``reduction of structure"
Theorem 1.1 from \cite{CGR}. The case where $\Char(k)$ divides
the order of the Weyl group of $G$ is particularly delicate;
here we use a refined version of \cite[Theorem 1.1]{CGR},
which is proved in~\cite{CGR2}.

The following symbols will be used for the remainder of the paper.
\[ \begin{array}{lcl}
 k          & & \text{algebraically closed base field of characteristic
                      $p \ge 0$} \\
 G      & & \text{connected reductive linear group defined
                                                       over $k$} \\
\overline{G} & & \text{a regular compactification of $G$} \\
 H       & & \text{finite abelian group} \\
 \ed       & & \text{essential dimension over k} \\
F_n = k((t_1)) \dots ((t_n)) & & \text{iterated Laurent series
field in $n$ variables}
\end{array} \]

\section{Regular compactifications}
\label{sect.wonderful}

Let $G$ is a connected reductive algebraic group defined over $k$.
Let $B$ and $B^-$ be opposite Borel subgroups of $G$.
By the Bruhat decomposition,
$G$ has finitely many $B \times B^-$ orbits. Hence, according to
\cite{brion-kumar} (Proposition 6.5), $G$
(viewed as a $G \times G$-variety) has a ``regular"
compactification, in the sense of \cite{bdp}; we will denote
this compactification by $\overline{G}$.
(Note that
the terms ``regular" and ``smooth" are not interchangeable
in this context; a regular compactification is smooth but
not the other way around.)
In particular, if $G$ is adjoint then $\overline{G}$ is the
wonderful compactification of $G$ constructed in~\cite{dp}
(and in prime characteristic in~\cite{strickland}).

Regular compactifications have many interesting special
properties. Most of them will not be used in the sequel.
The only property of $\overline{G}$ we will need
is the following description of the stabilizers of points in
$\overline{G}$ from~\cite[p. 151]{brion1}.

To every $\overline{g} \in \overline{G}$ we associate a pair of
opposite parabolic subgroups $P$ and $Q$, where $P$ is the projection of
$\Stab_{G \times G}(\overline{g})$ to the first factor of $G \times G$
and $Q$ is the projection to the second factor. Denote the unipotent
radicals of $P$ and $Q$ by $P_u$ and $Q_u$ and their common
Levi subgroup $P \cap Q$ by $L$.
The stabilizer $\Stab_{G \times G}(\overline{g})$ is then
given by
\begin{equation} \label{e.stabilizer}
\Stab_{G \times G}(\overline{g}) = \{ (p_u l z, q_u l)
\; | \; p_u \in P_u , \; q_u \in Q_u , \; l \in L, \;
z \in Z(L)^0 \} \, . \end{equation}

\begin{prop} \label{prop.abelian}
Let $F$ be a finite subgroup of $\Stab_{G \times G}(\overline{g})$
whose order is prime to $\Char(k)$. Then
$\mid\,  \rank( \pi_2(F)) -  \rank( \pi_1(F)) \, \mid \, \, \leq \rank(Z(L)^0)$.
\end{prop}

Here $\pi_1 \colon P \times Q \to P$ denotes the projection to the first
factor, $\pi_2 \colon P \times Q \to Q$ the projection to the second factor,
and $\rank(F)$ denotes the maximal value of $\rank(A)$, as $A$ ranges over
the abelian subgroups of a finite group $F$.

\begin{proof} The proof is based on tracing $F$ through the diagram
of natural projections
\[ \xymatrix{  &  \Stab_{G \times G}(\overline{g})
 \ar@{->}[dl]_{p_1} \ar@{->}[dr]^{p_2} &  \cr
P \ar@{->}[dr]_{\alpha_P}  &  &
Q \ar@{->}[dl]^{\alpha_Q}  \cr
 & L \ar@{->}[d] &  \cr
 & L/Z^0(L) , &  } \]
where $\alpha_P(p_u l) = l$ and $\alpha_Q(q_u l) = l$ for any
$p_u \in P_u$, $q_u \in Q_u$ and $l \in L$. Since
the kernels $P_u$ and $Q_u$ of $\alpha_P$ and $\alpha_Q$
are unipotent, we see
that $\alpha_P$ and $\alpha_Q$ project the finite groups
$\pi_1(F)$ and $\pi_2(F)$ isomorphically onto subgroups of $L$,
which we will denote by $F_1$ and $F_2$, respectively.
By~\eqref{e.stabilizer}, $F_1$ and $F_2$ have the same
image in $L/Z(L)^0$, which we will denote by $F_0$.
Since the natural projection
$\phi_{|F_2} \colon F_2 \to F_0$ is surjective
with kernel  $Z_2 \subset Z^0(L)$, we have
\[ \rank(F_2)  \le  \rank(F_1) +  \rank(Z_2)
\leq \rank(F_1) + \rank \, Z^0(L) \,.   \]
By symmetry,  we also have the reverse inequality
$ \rank(F_1) - \rank(F_2) \le \rank \, Z^0(L)$, and
the proposition follows.
\end{proof}

\section{Compactifications of homogeneous spaces}

The following lemma is well known in characteristic zero
(see \cite[Lemma 2.1]{ry2}).

\begin{lem} \label{lem.geometric-quotient}
Let $F$ be a finite group and
let $X$ be a normal quasiprojective $F$-variety. Then

\smallskip

(a) $X$ is covered by affine open $F$-invariants subsets.

\smallskip

\noindent Assume moreover that the  order of $F$ is invertible in $k$. Then

\smallskip

(b) there is a geometric quotient map $\pi: X \to X / F$.

\smallskip

(c) Moreover, if $X$ is projective, then so is $X / F$.

\end{lem}

Recall that we are assuming throughout that the base field $k$ is
algebraically closed.

\begin{proof}

\smallskip

\noindent (a) The proof of {\it loc. cit} is characteristic free.

\smallskip

\noindent (b) Recall first that the  group $F$ is linearly reductive,
see \cite[\S 1]{GIT}.  If $X$ is affine, part (b) is proved
in~\cite[Theorem 1.1 and Amplification 1.3]{GIT}.
The general case follows from part (a) and the
characteristic free glueing assertion in~\cite[Theorems 4.14]{pv}.

\smallskip

\noindent (c) See \cite[Theorem 3.14]{N}.
\end{proof}

Let $G/k$ be a connected reductive linear algebraic group,
$\overline{G}$ be a regular compactification, and
$F$ be a finite subgroup of $G$ of order prime to $\Char(k)$.
 We shall denote by
$\overline{G}/F$ the geometric quotient
of $\overline G$ for the action of the finite group $F$
(on the right).  It may be viewed as a (possibly
singular) compactification of the homogeneous space
$G/F$. By the properties of geometric quotients

\smallskip
(i) the fibers of the natural projection $\overline{G} \to \overline{G}/F$
are the $F$-orbits of the right action of $F$ on $\overline{G}$, and

\smallskip
(ii) the left action of $G$ on $\overline{G}$ descends
to $\overline{G}/F$.

\begin{prop} \label{prop.homog}
Let $G$ be a connected reductive group,
$\overline{G}$ be a regular compactification and $H_1, H_2$ be finite
abelian subgroups of $G$ whose orders are prime to $\Char(k)$.
 If $\overline{G}/H_2$ has an $H_1$-fixed
point then
\[ \rank(H_1) - \rank(H_2) \le \rank \, C_G(H_1)^0  \, . \]
\end{prop}

\begin{proof} Let $y$ be an $H_1$-fixed point of $\overline{G}/H_2$,
and $x$ be a point of $\overline{G}$ lying above $y$.  Let
$H = \Stab_{H_1 \times H_2}(x)$.  The fact that $y$ is
$H_1$-fixed means that for every $h_1 \in H_1$, there is
an $h_2 \in H_2$ such that $(h_1, h_2) \in \Stab_{G \times G}(x)$.
In other words, if we denote
by $\pi_i$ the natural projection $G \times G \to G$ to the $i$th factor
(where $i = 1$ or $2$) then $\pi_1(H) = H_1$.

Now let $P = \pi_1(\Stab_{G \times G}(x))$, $Q =
\pi_2(\Stab_{G \times G}(x))$ be a pair of opposite parabolics and
$L = P \cap Q$ be their common Levi subgroup, as in the previous section.
By Proposition \ref{prop.abelian} applied to the group $H$,  we have
\[ \rank(H_1) - \rank(\pi_2(H))    \le \rank \, Z(L)^0 \, . \]
Since $\rank(H_2) \ge \rank(\pi_2(H))$, we have
\[ \rank(H_1) - \rank(H_2)    \le \rank \, Z(L)^0 \, . \]
It remains to show that
\[ \rank \, Z(L)^0 \le \rank \, C_G(H_1)^0 \, . \]
By the Levi decomposition, $H_1 \subset P$ is conjugate (in $P$
and hence, in $G$) to a finite subgroup of $L$. Denote this subgroup
by $H_L$.  The connected centralizers $C_G(H_1)^0$
and $C_G(H_L)^0$ will then also be conjugate in $G$, and since
$Z^0(L) \subset C_G(H_L)^0$, we see that
\[ \rank \, C_G(H_1)^0 = \rank \, C_G(H_L)^0 \ge \rank \, Z^0(L) \, , \]
as claimed.
\end{proof}

\section{Proof of Theorem~\ref{thm.main} in characteristic zero}
\label{sect.char0}

The following lemma is well known; we supply a short proof
for lack of a direct reference.

\begin{lem} \label{lem.fixed-pts}
Consider a faithful action of a finite abelian group $A$
on an algebraic variety $X$, defined over a field $k$.
Assume $\Char(k)$ does not divide $|A|$. If $A$ fixes
a smooth $k$-point in $X$ then $\dim(X) \ge \rank(A)$.
\end{lem}

\begin{proof} Let $x \in X(k)$ be a smooth $A$-fixed $k$-point. Then $A$
acts on the regular local ring $R = \mathcal{O}_x(X)$ and on
its maximal ideal $M = \mathcal{M}_{x}(X)$.

Assume the contrary; $\dim(X) = d$ and $\rank(A) > d$. Then
the $A$-representation on the $d$-dimensional contangent space
$T_x(X)^* = M/M^2$ cannot be faithful; denote its kernel by
$A_0 \ne \{ 1 \}$.
Since $|A|$ is prime to $\Char(k)$, the map $M \to M/M^2$ of
$A$-representations splits. Thus $R$ has a system
of local parameters $t_1, \dots, t_d \in M$ such that each $t_i$
is fixed by $A_0$. Then $A_0$ acts trivially on the completion
$\widehat{R} = k[[t_1, \dots, t_d]]$, hence, on
$R \subset \widehat{R}$, and consequently on $X$.
This contradicts our assumption that the $A$-action on $X$
is faithful.
\end{proof}

For the remainder of this section we will assume that $k$ is
an algebraically closed field of characteristic zero.
Before proceeding with the proof of Theorem~\ref{thm.main}, we
recall that every $G$-variety
is birationally isomorphic to a smooth projective $G$-variety;
cf.~\cite[Proposition 2.2]{ry2}.  This fact (whose proof
relies on equivariant resolution of singularities and
thus requires the assumption that $\Char(k) = 0$)
will be used repeatedly for the remainder of this section.

\begin{lem} \label{lem1}
Every irreducible generically free $G$-variety is birationally
isomorphic to a projective $G$-variety $X$ of the form
$(\oG \times Y)/S$, where

\smallskip
(a) $Y$ is a primitive smooth projective $S$-variety,

\smallskip
(b) $S$ acts on $\overline{G} \times Y$ by $s \cdot (h, y) =
(hs^{-1}, s \cdot y)$, and
$(\overline{G} \times Y)/S$ is the geometric quotient for this action.
\end{lem}

Here by saying that $Y$ is a {\em primitive} $S$-variety, we mean
that $S$ transitively permutes the irreducible components of $Y$.
Note also that $S$ is a finite group acting on the smooth projective variety
$\overline{G} \times Y$. In this situation a geometric quotient
exists and is projective; see, e.g.,~\cite[Theorem 4.14]{pv} or
Lemma~\ref{lem.geometric-quotient}.

\begin{proof}
Let $G$ be a reductive group. By \cite[Theorem 1]{CGR} there exists
a finite subgroup $S \subset G$ such that $H^1(K, S) \to H^1(K, G)$
is surjective for every $K/k$.

Now recall that elements of $H^1(K, G)$ are in
a natural 1-1 correspondence with birational
isomorphism classes of a generically free
primitive $G$-varieties $X$, with $K = k(X)^G$
see, e.g., \cite[Section 1.3]{popov}. Thus $X$ is induced from
some generically free $S$-variety $Y$.  In other
words, $X$ is birationally isomorphic to
$(G \times Y)/S$.  Since the $S$-variety $Y$
is only defined up to birational isomorphism, we may assume
without loss of generality that it is smooth and projective.
Moreover,  $(G \times Y)/S$ is clearly
isomorphic to $(\oG \times Y)/S$, and the lemma follows.
\end{proof}

\begin{lem} \label{lem.dimension}
Let $G$ be a reductive group, $H \subset G$ be
a finite abelian subgroup, and $X$ be a generically free $G$-variety.
If $H$ fixes a smooth $k$-point of $X$ then
\[ \dim \, k(X)^G \ge \rank(H) - \rank \, C_G(H)^0 \, . \]
\end{lem}

\begin{proof}
By Lemma~\ref{lem1} there is a birational $G$-equivariant
isomorphism $X \stackrel{\simeq}{\dasharrow}(\oG \times Y)/S$,
where $Y$ is a smooth complete
$S$-variety and $S \subset G$ is a finite subgroup,
as in Lemma~\ref{lem1}. Note that
$\dim(Y) = \trdeg_k \, k(X)^G$; we will denote this number by $d$.
By the Going Down Theorem~\cite[Proposition A.2]{ry}
$(\oG \times Y)/S$ also has an $H$-fixed $k$-point; denote it by
$x = [\overline{g}, y]$. The fiber of the natural projection
$X = (\oG \times Y)/S \to Y/S$ containing $x$ is easily seen to be
$G$-equivariantly isomorphic to $\oG/S_0$, where $S_0 := \Stab_S(y)$.
Now observe that since $y$ is a smooth $k$-point of $Y$,
$S_0$ can contain no abelian subgroup of rank $\ge d$.
In other words, $\rank(S_0) \le d$. By Proposition~\ref{prop.homog},
\[ \rank(H) - \rank(S_0) \le \rank \, C_G(H)^0  \]
and thus
\[ d \ge \rank(S_0) \ge \rank(H) -  \rank \, C_G(H)^0 \, . \]
This completes the proof of Lemma~\ref{lem.dimension}.
\end{proof}

We are now ready to proceed with the proof of Theorem~\ref{thm.main}.
Let $V$ be a generically free linear $k$-representation of $G$.

\smallskip
(a) The essential dimension $\ed(G)$ is the minimal value of 
$\trdeg_k \, k(X)^G$, where the minimum is taken over all
dominant rational $G$-equivariant maps
\[ V \dasharrow X \, , \]
such that $X$ is a generically free $G$-variety; 
see~\cite[Section 3]{reichstein1}. Thus our goal is 
to show that \begin{equation} \label{e.ed}
\trdeg_k \, k(X)^G \ge \rank(H) - \rank \, C_G(H)^0 \, .
\end{equation}
After replacing $X$ by a birationally equivalent $G$-variety,
we may assume that $X$ is smooth and projective.
Since $V$ has a smooth $H$-fixed $k$-point (namely the origin),
the Going Down Theorem~\cite[Proposition A.2]{ry} tells us that
$X$ also has an $H$-fixed $k$-point (which is smooth, because
every $k$-point of $X$ is smooth). The inequality~\eqref{e.ed}
now follows from Lemma~\ref{lem.dimension}.

\smallskip
(b) It suffices to show that the inequality~\eqref{e.ed} holds if
there is a diagram of dominant rational $G$-equivariant maps of the form
\[ \xymatrix{ V' \ar@{-->}[d] \ar@{-->>}[dr] &   \cr
V & X \, , } \]
where $X$ is a generically free $G$-variety, $\dim(V') = \dim(V)$
and $[k(V'): k(V)]$ is prime to $l$. (Note that the $G$-variety $V'$ 
is not required to be linear.)
Once again, we may assume without loss of generality that $V'$ and
$X$ are smooth and complete.  Since $H$ fixes the origin in $V$,
the Going Up Theorem~\cite[Proposition A.4]{ry} tells us that
$V'$ has an $H$-fixed $k$-point. Now by the Going Down
Theorem~\cite[Proposition A.2]{ry}, $X$ has an $H$-fixed
$k$-point as well, and Lemma~\ref{lem.dimension} completes the proof.
\qed

\section{The field of iterated Laurent series}

In this section we will describe the structure of
the iterated Laurent polynomial field
$F_n = k((t_1)) \dots ((t_n))$. 
Our proof of Theorem~\ref{thm.main} in full generality (i.e., without
assuming that $\Char(k) = 0$) will make use of these results.

We begin by describing the absolute Galois group
$\Gal(F_n)$.
Since $F_n=F_{n-1}((t_n))$ is a complete valuated field
with residue field $F_{n-1}$, we have an  exact sequence
$$
1 \to I_n \to \Gal(F_n)     \buildrel \pi_n \over \to \Gal(F_{n-1}) \to 1
$$
where $I_n$ stands for the inertia group; see \cite[\S II.7]{GMS}.
Let $\widehat \bbZ'$ be the prime to $p$ part
of $\widehat \bbZ$, i.e $\widehat \bbZ' =
\prod_{q \not = p} \widehat \bbZ_q$. In particular, if
$p = 0$ then $\widehat \bbZ' = \widehat \bbZ$.

\begin{lem}\label{lem.inertia}  There is a split exact sequence
\begin{equation} \label{e.inertia}
1 \to J_n \to   \Gal(F_n)  \to  (\widehat \bbZ' )^n \to 1 .
\end{equation}
such that
\begin{enumerate}

\item   $J_n=1$ if $p=0$,

\item $J_n$ is a free pro-p-group if $p >0$.
\end{enumerate}
\end{lem}

\begin{proof} We proceed by induction on $n$.
The group $I_n$ fits in an exact sequence
$$ 0 \to I_n^{wild} \to I_n \to \widehat \bbZ' \to 1 \, , $$
where $I_n^{wild}$ is the wild inertia group (it is a pro-p-group).
Define
$F_{n,m}:= k(( \root m \; \; \of t_1)) \dots (( \root m \; \; \of t_n))$ and
$F_{n,\infty}:= \limind\limits_{(m,p)=1} F_{n,m}$.
Since \[ \Gal(F_{n,\infty}/ F_{n-1,\infty}((t_n))) =
\widehat \bbZ' \, , \]
we have the following commutative diagram of profinite groups
{\small
$$
\begin{CD}
&& 1 && 1 &&1 \\
&& @VVV   @VVV  @VVV \\
1 @>>> I_n^{wild}  @>>> J_n @>>> J_{n-1} @>>> 1 \\
&& @VVV   @VVV  @VVV \\
1 @>>> I_n  @>>> \Gal(F_n) @>>> \Gal(F_{n-1}) @>>> 1 \\
&& @VVV   @VVV  @VVV \\
 1@>>> \Gal(F_{n,\infty}/ F_{n-1,\infty}((t_n))) @>>>
\Gal(F_{n,\infty}/ F_{n}) @>>> \Gal(F_{n-1,\infty}/ F_{n-1}) @>>> 1 .\\
&& @VVV   @VVV  @VVV \\
&& 1 && 1 &&1 \\
\end{CD} $$ }

If $p=0$, it follows that $J_n=1$.
If $p>0$, we see by induction that  $J_n$ is pro-p-group.
The group $J_n$ is the absolute Galois group
of the field $F_{n,\infty}$, so ${\rm cd}_p(J_n) \leq 1$
(\cite{serre-gc}, II.3.1, proposition 7) and
$J_n$ is then a free pro-p-group ({\it ibid}, I.4.2, Corollary 2).

Since $(\widehat \bbZ')^n$ is a prime-to-$p$ group, we conclude that the
sequence~\eqref{e.inertia} splits; see~\cite[I.5.9, Corollary 1]{serre-gc}.
\end{proof}

We will now show that every finite field extension 
of $F_n$ is $k$-isomorphic to $F_n$.
Recall that the lexicographic order $\prec$ on  ${\bf Z}^n$ is
defined as follows:
\[ (m_1, \cdots m_n) \prec  (m'_1, \cdots m'_n) \]
if $m_i < m'_i$ for the smallest subscript $i$ with $m_i \not = m'_i$.
A valuation $v=(v^{(n)}, \cdots ,v^{(1)}): E^\times \to \Gamma$
on a field $E$ is called $n$-{\em discrete} if the group
$\Gamma$ is isomorphic to a (lexicographically ordered) subgroup
of ${\bf Z}^n$; see~\cite[1.1.3]{F}.
Then $E$ is a $1$-discrete valuation with respect to the
first component $v(^{n})$ of $v$ and the residue field
$E_{n-1}$ is a $(n-1)$-discrete valuation via $(v^{(n-1)}, \cdots ,v^{(1)})$.
In this way we obtain a sequence of fields
$E=E_n$, $E_{n-1}$, ... , $E_1$ such that
$E_i$ is the residue field of $E_{i+1}$ with respect
to a $1$-discrete valuation. The residue field $E_0$ of $E_1$ then
coincides with the residue field $\overline E_v$.

The definition of completeness for $E$ is inductive
as follows~\cite[1.2.1]{F}.
The field $E$ is a complete $n$-discrete field
if $E_n$ is complete to respect to $v^{(n)}$
and $E_{n-1}$ is complete.
Assume from now on that $E$ is a complete $n$-discrete field.
Then according to \cite[3.1]{wadsworth}, $E$  is henselian,
i.e its valuation ring is a henselian ring.
In particular, given a finite extension $E'/E$
 the valuation
  $v$ extends uniquely to $v' : {E'}^\times
\to \frac{1}{[E':E]} {\bf Z}^n$,  the formula being
$v' = \frac{1}{[E':E]}  \, v \circ N_{E'/E}$.
Then  $E'$ is a $n$-discrete field which is complete (by induction).

The field $F_n=k((t_1))((t_2)) \cdots ((t_n))$ over
iterated Laurent series over a base field $k$ is $n$-complete;
see \cite[\S 3]{wadsworth}. Here the valuation $v$ on $F_n$ is
defined by
$$
v\Bigl( \sum_{i_1} \dots \sum_{i_n} \,
c_{i_1, \dots,i_n} \, t_1^{i_1} \dots t_n^{i_n} \Bigr)
= {\rm Min}\Bigl\{ \, (i_1,...,i_n) \enskip
\mid c_{i_1,...,i_n} \not =0 \Bigr\}.
$$

\begin{prop}\label{structure}
Let $E$ be a complete $n$-discrete $k$-field. Then

\begin{enumerate}
\item $E$ is isomorphic to $\overline E_v((t_1))((t_2)) \cdots ((t_n))$;
\item Suppose the valuation $v$ is trivial on a perfect subfield $K$ of $E$.
Then $E$ is $K$--isomorphic to
$\overline E_v((t_1))((t_2)) \cdots ((t_n))$.
\end{enumerate}
\end{prop}

\begin{proof} (1) immediately follows from (2) if we
take $K$ to be the prime subfield of $E$.

To prove (2), first assume that $n=1$. In this case
Cohen's structure theorem \cite[Theorem 10]{Co} shows
that the valuation ring
of $E$ is $K$-isomorphic to $\overline E_v[[t_1]]$.
Thus $E$ is $K$-isomorphic to $\overline E_v((t_1))$.

Now suppose $n \ge 2$.
Let $E=E_n$, $E_{n-1}$,..., $E_1$ the sequence of fields
constructed above.  By induction, we may assume
that $E_{n-1}$ is $K$--isomorphic to
$\overline E_v((t_1))((t_2)) \cdots ((t_{n-1}))$.
Since we know part (2) holds for $n=1$, we conclude that
$E=E_n$ is $K$-isomorphic
to $\overline E_v((t_1))((t_2)) \cdots ((t_{n}))$.
\end{proof}

\begin{cor}\label{laurent} Assume that $k$ is algebraically
closed.  Then any finite extension of $F_n=k((t_1))...((t_n))$
is $k$--isomorphic to $F_n$.
\end{cor}

\begin{proof}
A finite extension $E$ of $F_n$ is a complete $n$-discrete field.
Its residue field $E_v$ is a finite extension of $k$, since $k$ is
algebraically closed, we conclude that $E_v = k$.
Proposition~\ref{structure}(2) now shows
that $E$ is $k$--isomorphic to $F_n$.
\end{proof}

\section{Reduction of structure}
\label{sect.reduction}

In this section $G/k$ will denote a linear algebraic group
defined over an algebraically closed field $k$ of characteristic $\geq 0$,
whose identity component $G^0$ is reductive.  $F_n$ will denote
the iterated power series field $k((t_1))((t_2))....((t_n))$
in variables $t_1, \dots, t_n$, as in the previous section.
As usual, we will say that $\gamma \in H^1(K, G)$

\smallskip
{\em descends} to a subfield $K_0 \subset K$ if $\gamma$ lies in the image
of the restriction map $H^1(K_0, G) \to H^1(K, G)$

\smallskip
{\em admits reduction of structure} to a subgroup $A \subset G$ if $\gamma$
lies in the image of the natural map $H^1(K, A) \to H^1(K, G)$.

\begin{prop}\label{prop.compress}
Suppose $\gamma \in H^1(F_n, G)$ descends to a subfield
$K \subset F_n$ such that $\trdeg_k(K) = d < \infty$.

\smallskip
(a) Assume that $\Char(k) = 0$.
Then $\gamma$ admits reduction of structure
to a finite abelian subgroup $A \subset G$ of rank $\le d$.

\smallskip
(b) Assume that $\Char(k) = p > 0$.
Then there exists a finite field extension $F'/F_n$ such that
$[F': F_n]$ is a power of $p$ and $\gamma_{F'}$ admits
reduction of structure to a finite abelian subgroup
$A \subset G$ of rank $\le d$, whose order $|A|$ is prime to $p$.
\end{prop}

Our proof of Proposition~\ref{prop.compress} will make use of the following
two simple lemmas.

\begin{lem} \label{lem.galois}
Suppose $K \subset E$ is a field extension
such that $K$ is algebraically closed in $E$.
Then

\smallskip
(a) for every finite Galois field extension $K'/K$,
$K'_E = K' \otimes_K E$ is a field.

\smallskip
(b) The absolute Galois group $\Gal(K)$ is
a quotient of $\Gal(E)$.
\end{lem}

\begin{proof}
(a) By the primitive element theorem we can write $K'$ as
$K[x]/(p(x))$ for some irreducible monic polynomial
$p(x) \in K[x]$.  Then $K'_E = K' \otimes_K E = E[x]/(p(x))$,
and we need to show that $p(x)$ remains irreducible over $E$.

We argue by contradiction. Suppose $p(x) = p_1(x) p_2(x)$
for some non-constant monic polynomials $p_1(x), p_2(x) \in E[x]$.
The coefficients of $p_i(x)$ are then polynomials
in the roots of $p(x)$. (Here $i = 1$ or $2$.)
In particular, they are
algebraic over $K$. Since $K$ is algebraically
closed in $E$, we conclude that $p_i(x) \in K[x]$.
Thus $p(x)$ is reducible over $K$, a contradiction.

\smallskip
(b) Let $\overline{K}$ be the algebraic closure of $K$.  Then
$\Gal(K) = \Gal(\overline{K}/K) = \Gal(\overline{K}_E/E)$, and
$\overline{K}_E = \overline{K} \otimes_K E$
is an $E$-subfield of $\overline{E}$ by part (a).
\end{proof}

\begin{lem} \label{lem.gabber}
Let $d\geq 0$ be an integer.  Let $\Gamma$ be a finitely
generated abelian profinite group such that $\cd(\Gamma) \leq d$.
Then $\Gamma$ is a direct summand of ${\widehat \bbZ}^d$.
\end{lem}

\begin{proof} Without loss of generality, we may assume
that $\Gamma$ is a $p$-profinite group for a prime $p$.
Since $\Gamma$ is finitely generated and abelian,
$\Gamma \cong \bbZ_p^m$ for some integer $m$;
see \cite[Theorem 4.3.4.(a)]{RZ}.
Thus $\cd(\Gamma)=m$. Since we are assuming $r \leq d$, this shows that
$\Gamma$ is a direct summand of ${\widehat \bbZ}^d$.
\end{proof}

\begin{proof}[Proof of Proposition~\ref{prop.compress}]
Let $W$ be the Weyl group of $G$.
We recall that there exists a finite subgroup
$S \hookrightarrow G$ such that every prime factor of $|S|$
divides $|W|$ and $S$ has the following property.

\smallskip
(i) If ${\rm Char}(k)$ does not divide  $|W|$, then
the map $H^1(K,S) \to H^1(K,G)$ is surjective for every field $K/k$;
see \cite{CGR}.

\smallskip
(ii) If $\Char(k)$ divides $|W|$, the above map is surjective
for every perfect field $K/k$; see \cite{CGR2}.

\smallskip
We fix the finite subgroup $S$ with these properties for the rest of
the proof.

Let $\rho: (\widehat \bbZ' )^n  \to  \Gal(F_n)$ be
a splitting of the exact sequence~\eqref{e.inertia} in
Lemma~\ref{lem.inertia}.
Denote the extension of $F_n$ associated by
the Galois correspondence to the image of $\rho$
by $E_n/F_n$ and its perfect closure by $E_n^{perf}/F_n$.
Note that if $\Char(k) = 0$ then $\rho$ is an isomorphism and
$E_n^{perf} = E_n = F_n$.
If $\Char(k) = p$ then the degree of any finite subextension
of $E_n/F_n$ is a power of $p$ and $\Gal(E_n) = (\widehat \bbZ' )^n$.
The same is true for the perfect closure $E_n^{perf}/F_n$.

Since Galois cohomology commutes with direct limits of fields,
in order to establish parts (a) and (b) of the proposition,
it suffices to show
that $\gamma_{E_n ^{perf}}$ admits reduction of structure
to some abelian subgroup $A \subset S$ of rank $\le \trdeg_k(K)$
(where $|A|$ is prime to $p$, if $p = \Char(k) > 0$).
After replacing $K$ by its algebraic closure in $E_n^{perf}$,
we may assume that $K$ is algebraically closed
in $E_n^{perf}$.  In particular, $K$ is perfect.

By our assumption $\gamma$ descends to some
$\gamma_K \in H^1(K, G)$. On the other hand, by (ii)
$\gamma_K$ is the image of some $\delta_K \in H^1(K, S)$.
The class $\delta_K$ is represented by a continuous
homomorphism $\psi \colon \Gal(K) \to S$. Clearly $\delta_K$
(and hence, $\gamma_K$ and $\beta_{E_n^{perf}}$)
admit reduction of structure to the subgroup
$A = \Ima(\psi)$ of $S$.  It remains to show that $A$ is
an abelian group of rank $\le d$ whose order is prime to $p$.

By Lemma \ref{lem.galois} we can identify $\Gal(K)$ with
a quotient of $\Gal(E_n^{perf}) = (\widehat \bbZ' )^n$.
In particular, $\Gal(K)$ is finitely generated,
abelian, and the order of every finite quotient of $\Gal(K)$
is prime to $p$.  Moreover, by Tsen's theorem,
$\cd(\Gal(K)) \le d$; cf.~\cite[II.4.2]{serre-gc}.
Thus  Lemma~\ref{lem.gabber} enables us to conclude that
$\Gal(K)$ is a direct summand of $(\widehat{\bbZ}')^d$.
Hence, the finite quotient $A$ of $\Gal(K)$ is an abelian
group of rank $\le d$ whose order is prime to $p$.
\end{proof}

\begin{remark} \label{rem.compress}
A minor modification of the above argument
(in particular, using (i) instead of (ii)) shows that
the assertion of Proposition~\ref{prop.compress}(a)
holds whenever ${\rm Char}(k)$
does not divide the order of the Weyl group $W$ of $G$.
In other words, in this case we can take $F'$ to be $F_n$
in part (b).  Since we will not use this result in
the sequel, we leave the details of its proof 
as an exercise for an interested reader.
\end{remark}

\section{Fixed points in homogeneous spaces}
\label{sect.fixed}

Let $k$ is an algebraically closed field of characteristic $p \ge 0$,
$t_1, \dots t_n$ are independent variables over $k$,
and $H = (\bbZ/ m \bbZ)^n$. If $p > 0$, we will assume that
$m$ is prime to $p$.
We we will continue to denote the iterated power series field
$k((t_1))((t_2)) \dots ((t_n))$ by $F_n$.

The purpose of this section is to establish the following fixed
point result which will be used in the proof of Theorem~\ref{thm.main}.
For notational convenience, we will consider an arbitrary
(not necessarily injective) morphism
\begin{equation} \label{e.phi}
\phi \colon H = (\bbZ/m\bbZ)^n \to G
\end{equation}
of algebraic groups.  This is slightly more general than
considering a finite abelian subgroup of $G$. We will assume
that $G$, $H$ and $\phi$ are fixed throughout this section.

\begin{prop} \label{prop.GP} Assume that

\smallskip
(1) $F'/F_n$ is a finite field extension of degree prime to $|H|$,

(2) $\beta \in H^1(F_n, H)$ is represented by an $H$-Galois field
extension $E/F$, and

\smallskip
(3) $\phi_*(\beta)_{F'} \in H^1(F', G)$ admits reduction of structure to
a finite subgroup $S$ of $G$.

\smallskip
Then $\phi(H)$ has a fixed $k$-point in any $G$-equivariant
compactification $Y$ of $G/S$.
\end{prop}

Here by a $G$-equivariant compactification of $G/S$ we mean
a complete (but not necessarily smooth) $G$-variety, which contains
$G/S$ as a dense open $G$-subvariety.

\begin{proof} By Corollary~\ref{laurent}, $F'$ is $k$-isomorphic
to $F_n$. Thus, after replacing $F_n$ by $F'$ and $\beta$
by $\beta_{F'}$, we may assume that $F' = F_n$. (Note that
$\beta_{F'} \in H^1(F', H)$ is represented by the $H$-Galois algebra
$E_{F'}/F'$, where $E_{F'} = E \otimes_{F} F'$. Since $E$ is a field
and $[F_n':F_n]$ is prime to $|H| = [E:F]$, $E_{F'}$ is again a field.)

By Lemma~\ref{lem.inertia}, we may assume that
$E= k((s_1)) \dots ((s_n))$, where $s_i^m = t_i$ and
there exists a minimal set of generators $\tau_1,..., \tau_n$
of $H$ such that $H$ acts on $k((s_1)) \dots ((s_n))$ by
  \begin{equation} \label{e.kummer}
  \tau_i(s_j)= \begin{cases} \text{$\zeta \, s_j$, if $i = j$ and} \\
  \text{$s_j$, if $i \neq j$, } \end{cases}
  \end{equation}
where $\zeta$ is a primitive $m$th root of unity (independent of $i$ and $j$).
In the sequel we will denote $E$ by $F_{n, m}$; note that we previously
encountered this field in the proof of Lemma~\ref{lem.inertia}.

Set $\gamma = \phi_*(\beta) \in H^1(F_n, G)$ and
consider the twisted $F_n$-variety ${_\gamma} Y$
which is a compactification of the twisted variety ${_\gamma} (G/S)$.
By our asumption ${_\gamma} (G/S)$ has a $F_n$-point, so a fortiori
$$
{_\gamma Y} (F_{n, m})= \Bigl\{ y \in Y(F_{n, m}) \, \mid \,
\gamma (\sigma). ^\sigma y =y \enskip \forall \sigma \in H
\Bigr\} \neq \emptyset .
$$
Since $Y$ is complete, this implies
${_\gamma} Y(F_{n-1,m}[[{\root m \of t_n}]]) \neq \emptyset$.
Specializing $t_n$ to $0$, we see that
$$
\Bigl\{ y \in Y(F_{n-1,m}) \, \mid \, \gamma(\sigma). ^\sigma
y =y \enskip \forall \sigma \in H\Bigr\} \not = \emptyset,
$$
where the Galois action of $H$ on
$Y(F_{n-1,m})$ is induced by
the canonical projection $H\to \Gal( F_{n,m}/F_{n})
\to \Gal( F_{n-1,m}/F_{n-1})$.
Repeating this process, we finally obtain
$$
\Bigl\{ y \in Y(k) \, \mid \, \gamma(\sigma). ^\sigma y =y
\enskip \forall \sigma \in H \Bigr\} \not = \emptyset.
$$
Since $k$ is algebraically closed, we conclude that $\phi(H)$
fixes some $k$-point of $Y$.
\end{proof}

\begin{cor} \label{cor.class} Let $k$ be an algebraically closed field
and $G/k$ be a connected reductive group.
Suppose there exists a class $\beta \in H^1(F_n, H)$ such that
$\beta$ is represented by an $H$-Galois field extension of $F_n$. If
$\phi_*(\beta) \in H^1(F_n, G)$ descends to a $k$-subfield
$K \subset F_n$ then
\[ \trdeg_k(K) \ge \rank \, \phi(H) - \rank \, C_G(\phi(H))^0 \, . \]
\end{cor}

Here $H$ and $\phi$ are as in~\eqref{e.phi}.

\begin{proof} Let $\trdeg_k(K) = d$.
By Proposition \ref{prop.compress} there exists
a finite extension $F'/F_n$ and a
finite abelian subgroup $A \subset G$ of rank $\le d$
such that

\smallskip
$|A|$ is prime to $\Char(k)$,

\smallskip
$F' = F_n$ if $\Char(k) = 0$, and $[F':F_n]$ is a power of $p$ if
$\Char(k) = p$,

\smallskip
and

\smallskip
$\phi_*(\beta)_{F'}$ admits reduction of structure to $A$.

\smallskip
Let $\overline{G}$ be a regular compactification of $G$.
By Proposition~\ref{prop.GP} $Y = \overline G/A$ has a $\phi(H)$-fixed point.
Now Proposition~\ref{prop.homog}, with $H_1 = \phi(H)$ and $H_2 = A$,
tells us that
\[ \rank \, \phi(H) - \rank \, A \le  \rank \, C_G(\phi(H))^0 \, . \]
Consequently,
\[ d \ge \rank \, A \ge \rank \,\phi(H) - \rank \, C_G(\phi(H))^0 \, , \]
as claimed.
\end{proof}

\section{Proof of Theorem~\ref{thm.main}}
\label{char.free(a)}

In the statement of Theorem~\ref{thm.main}, we assume that
$H$ is a subgroup of $G$, where as in the previous section we 
worked with a homomorphism $\phi \colon H \to G$ instead.
For notational consistency, we will restate Theorem~\ref{thm.main}
in the following (clearly equivalent) form.

\begin{thm} \label{thm.main2}
Let $G$ be a connected reductive linear algebraic group defined
over an algebraically closed base field $k$, $H \simeq (\bbZ/m \bbZ)^n$
and $\phi \colon H \to G$ be a (not necessarily injective) homomorphism
of algebraic groups. Assume
$\Char(k)$ does not divide $m$. Then

\smallskip
(a) $\ed(G) \ge \rank \, \phi(H) - \rank \, C_G(\phi(H))^0$.

\smallskip
(b) Moreover, if $H$ is a $l$-group (i.e., $m$ is a power of
a prime integer $l$)
then $\ed(G; l) \ge \rank \, \phi(H) - \rank \, C_G(\phi(H))^0$.
\end{thm}

Let $t_1, \dots, t_n$ be independent variables over $k$,
$K_n = k(t_1, \dots, t_n)$ and $K_{n, m} = k(s_1, \dots, s_n)$,
where $s_i^m = t_i$. The $H$-Galois field extension $K_{n, m}/K_n$
gives rise to a class $\alpha \in H^1(K_n, H)$. We will be interested in
the class $\phi_*(\alpha) \in H^1(K_n, G)$, which we will sometimes
refer to as a {\em loop torsor}.  (Such torsors naturally come up in
connection with loop algebras; see~\cite{GP}.) We are now ready
to proceed with the proof of Theorem~\ref{thm.main2}(a).

\smallskip
(a) By the definition of $\ed(G)$, it suffices to show that
\[ \ed \, \phi_*(\alpha) \ge \rank \,\phi(H) - \rank \, C_G(\phi(H))^0 \, . \]
Let $d = \ed(\phi_*(\alpha))$. Let $\beta = \alpha_{F_n} \in H^1(F_n, H)$.
Then $\beta$ is represented by the field extension $F_{n, m}/F_n$.
Moreover, $\phi_*(\beta)$ descends to $\phi^*(\alpha) \in H^1(K_n, G)$,
which, by our assumption, further descends to a $k$-subfield of $K_n$
of transcendence degree $d$. Corollary~\ref{cor.class} now tells us that
\[ d \ge \rank \, \phi(H) - \rank \, C_G(\phi(H))^0 \, . \]
This completes the proof of Theorem~\ref{thm.main2}(a).

Our proof of Theorem~\ref{thm.main2}(b) will rely
on the following lemma.

\begin{lem} \label{lem.primetol}
Let $l$ be a prime integer, $E/K$ be a finite
field extension of degree prime to $l$ and $F/K$ be an arbitrary
(not necessarily finitely generated) field extension.
Assume $\Char(K) \ne l$. Then there exists a field extension $F'/F$
and a $K$-embedding $E \hookrightarrow F'$ such that the diagram
\[ \xymatrix{E \; \ar@{-}[dd]_{\text{\tiny prime-to-$l$}} \ar@{^{(}->}[r] &
                   F' \ar@{-}[d]^{\text{\tiny prime-to-$l$}} \cr
                      & F \cr
               K \ar@{-}[ur]   &   } \]
commutes.
\end{lem}

\begin{proof} We may assume without loss of generality that $E=K(\alpha)$.
Let $P$ be the minimal polynomial of $\alpha$
and let $P=  Q_1^{e_1} ... Q_r^{e_r}$ be the prime  decomposition of
$P$ as $F$-polynomial.
We have
$$
F \otimes_K E = F \otimes_K { {K[X]} \over {P}}
\cong   { {F[X]} \over {Q_1^{e_1}}} \oplus \cdots \oplus { {F[X]} \over {Q_r^{e_r}}}.
$$
Denote by $F_i=F[X]/Q_i(X)$ the residue field of the prime polynomial $Q_i$.
  Since
\[ [K':K]=  [F_1:F]^{e_1} + \dots + [F_r:F]^{e_r} \]
is prime to $l$, one of the degrees $[F_i : F]$ is prime to $l$
for some $i = 1, \dots, r$. We can now take $F' = F_i$.
\end{proof}

We are now ready to finish the proof of Theorem~\ref{thm.main2}(b).
Once again, we will denote
the class of the $H$-Galois field extension $K_{n, m}/K_n$ in $H^1(K_n, H)$
by $\alpha$ and consider  the loop torsor $\phi_*(\alpha) \in H^1(K_n, G)$.
By the definition of $\ed(G; l)$ it suffices to show that
\[ \ed(\phi_*(\alpha); l) \ge \rank \, H - \rank \, C_G(H)^0 \, . \]
Equivalently, we want to show that
\[ \ed(\phi_*(\alpha)_E) \ge \rank \, H - \rank \, C_G(H)^0 \]
for every finite extension $E/K_n$ of degree prime to $l$.
Suppose $E/K_n$ is such an extension and $\ed(\phi_*(\alpha)_E) = d$.

By Lemma~\ref{lem.primetol} there exists a finite field extension
$F'/F_n$ of degree prime to $l$ and
a $K_n$-embedding $\colon E \hookrightarrow F'$ such that
the diagram
\[ \xymatrix{E \; \ar@{-}[dd]_{\text{\tiny prime-to-$l$}} \ar@{^{(}->}[r] &
                   F' \ar@{-}[d]^{\text{\tiny prime-to-$l$}} \cr
                      & F_n \cr
               K_n \ar@{-}[ur]   &   } \]
commutes.  We want to conclude that
$d \ge \rank \, H - \rank \, C_G(H)^0$
by applying Corollary~\ref{cor.class} to
$\beta = \alpha_{F'} \in H^1(F', H)$.
Since $F'$ is $k$-isomorphic to $F_n$ (see
Corollary~\ref{laurent}), Corollary~\ref{cor.class}
can be applied to this $\beta$, as long as we can show that

\smallskip
(i) $\phi_*(\beta)$ descends to a $k$-subfield of $F'$, of
transcendence degree $d$ over $k$, and

\smallskip
(ii) $\beta$ is represented by an $H$-Galois field extension of $F'$.

\smallskip
\noindent
(i) is clear since $\phi^*(\beta) \in H^1(F', G)$ descends
to $\phi^*(\alpha) \in H^1(E, G)$, which, by our assumption, descends to
a $k$-subfield $E_0 \subset E$ with $\trdeg_k(E_0) = d$. To prove (ii),
note that $\alpha_{F_n}$ is represented by the field extension $F_{n, m}/F_n$.
Thus $\beta = \alpha_{F'}$ is represented by the $H$-Galois algebra
$F_{n, m} \otimes_{F_n} F'$ over $F'$.
Since $[F':F_n]$ is a finite and prime to $l$,
$F_{n, m} \otimes_{F_n} F'$ is a field.
This concludes the proof of (ii) and thus of Theorem~\ref{thm.main2}.
\qed

\section{Examples}
\label{sect.examples}

This section contains five examples illustrating Theorem~\ref{thm.main}(b).

\begin{example} \label{ex}
If $\Char(k) \ne 2$ then $\ed(\GO_n; 2) \ge n - 1$.
\end{example}

\begin{proof}
Let $H \simeq (\bbZ/2 \bbZ)^n$ be the subgroup of diagonal
matrices in $\Orth_n$. Viewing $H$ as a subgroup of $\GO_n$,
we easily see that $C_{\GO_n}^0(H)$  = the center of $\GO_n$,
has rank $1$.  Applying Theorem~\ref{thm.main}(b) to this subgroup
we obtain the desired bound.
\end{proof}

\begin{example} \label{ex2}
If $p$ is a prime and $\Char(k) \neq p$ then

\smallskip
(a) $\ed(SL_{p^r}/\mu_{p^s}; p) \ge \begin{cases}
\text{$2s + 1$, if $s < r$,} \\
\text{$2s$, if $s = r$} \end{cases}$ and

\smallskip
(b) $\ed(GL_{p^r}/\mu_{p^s}; p) \ge \begin{cases}
\text{$2s$, if $s < r$,} \\
\text{$2s - 1$, if $s = r$.} \end{cases}$
\end{example}

\begin{proof} (a) The group $\SL_{p^r}/\mu_s$ has a self-centralizing
subgroup $H \simeq (\bbZ/p \bbZ)^r \times \bbZ/p^{r - s} \bbZ$; see
\cite[Lemma 8.12]{ry}. Now apply Theorem~\ref{thm.main}(b) to this group.

\smallskip
(b) We now view $H$ as a subgroup of
$\GL_{p^r}/\mu_{p^s}$. The centralizer
$C_{\GL_{p^r}/\mu_s}^0(H)$ is  the center of $\GL_{p^r}/\mu_s$;
it is isomorphic to a 1-dimensional torus.
Part (b) now follows from Theorem~\ref{thm.main}(b).
\end{proof}

The non-vanishing of the Rost invariant
$H^1(\ast \,, G) \to H^3(\ast \, , \mu_p)$ for a group $G$
and a prime $p$
implies that $\ed(G; p) \ge 3$; cf.~\cite[Theorem 12.14]{reichstein1}.
In particular, one can show that $\ed(F_4; 3), \ed(E_6; 2) $ and
$\ed(E_7; 3) \ge 3$ in this way.
In Examples~\ref{ex.F_4}--\ref{ex.E_7} below we will deduce
these inequalities directly from Theorem~\ref{thm.main}(b)
and show that equality holds in each case. 

\begin{example} \label{ex.F_4}
If $\Char(k) \ne 3$ then $\ed(F_4; 3) = 3$.
\end{example}

\begin{proof}
$F_4$ has a self-centralizing subgroup isomorphic to $(\mu_3)^3$;
see~\cite[Theorem 7.4]{griess}. Theorem~\ref{thm.main} now tells us that
$\ed(F_4;3) \geq 3$. To prove the opposite inequality, recall that
$H^1(K, E_6)$ classifies
the exceptional 27-dimensional Jordan algebras
$J/K$.  After a quadratic extension $K'/K$,
${J \otimes_K K'}$ is given by the first Tits construction.
Without loss of generality, we may assume that $J=(A, \nu)$
where $A$ is a central simple $K$-algebra of degree $3$ and $\nu$
is a scalar in $K$. Since $A$ is a symbol algebra $(a, b)_3$, we conclude that
that ${J \otimes_K K'}$ descends to the subfield $k(a, b, \nu)$ of $K$ of
transcendence degree $\le 3$. We conclude that $\ed(J;3) \leq 3$ and
thus $\ed(F_4; 3) \le 3$, as claimed.
\end{proof}

\begin{example} \label{ex.E_6}
If $\Char(k) \ne 2$ then $\ed(E_6; 2) = 3$.
\end{example}

Here $E_6$ denotes the simply connected simple groups of type $E_6$.
By abuse of notation we will also write $E_6$ for the Dynkin diagram of $E_6$.

\begin{proof}
By \cite[Table II]{griess}, $E_6$ has a unique non-toral
subgroup $H$ isomorphic to $(\bbZ/ 2 \bbZ)^5$.
To compute the rank of its centralizer, we make use
of its Witt--Tits index $I(H) \subset E_6$ which
describes the type of a minimal parabolic subgroup
containing $H$; see \cite[Section 3]{GP}.
The Dynkin diagram for $E_6$ is

\smallskip

$$
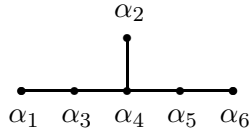
\begin{array}{ll}
\begin{picture}(85,30)
\put(00,02){\line(1,0){80}}
\put(40,02){\line(0,1){20}}
\put(00,02){\circle*{3}}
\put(20,02){\circle*{3}}
\put(40,02){\circle*{3}}
\put(60,02){\circle*{3}}
\put(80,02){\circle*{3}}
\put(40,22){\circle*{3}}
\put(00,02){\circle*{3}}
\put(80,02){\circle*{3}}
\put(-5,-10){$\alpha_1$}
\put(15,-10){$\alpha_3$}
\put(35,-10){$\alpha_4$}
\put(55,-10){$\alpha_5$}
\put(75,-10){$\alpha_6$}
\put(35,30){$\alpha_2$}
\end{picture}
\end{array}
$$

\bigskip

\smallskip

\noindent Set $I= \{ \alpha_2, \alpha_3, \alpha_4, \alpha_5 \}$. Let
$P_I$ be a standard parabolic subgroup and
$L_I = Z_G(T_I)$ its standard Levi subgroup.  Then $DL_I=  \Spin_8$.
Since  $\Spin_8$ has a  maximal non-toral  2-elementary abelian subgroup
of rank 5 (see \cite[Table I]{griess}), we may assume that $H \subset \Spin_8$.
Moreover, $C_{\Spin_8}(H)$ is finite, so $H$ is irreducible in $L_I$.
It follows that $P_I$ is a minimal parabolic subgroup
of $E_6$ containing H; the Witt--Tits index of $H$ is then $I$.
By \cite[Proposition 3.19]{GP}, we have
\[ \rank \, C_{E_6}(H) = | E_6 \setminus I |  = 2 \, . \]
Theorem~\ref{thm.main}(b) now tells us that $\ed(E_6; 2) \ge 5-2= 3$.

To prove the opposite inequality, suppose
$\alpha \in H^1(K,E_6)$, where $K$ is a field containing $k$.
Let $L \subset E_6$ be the Levi subgroup of
the parabolic $E_6 \setminus \{\alpha_6 \}$.
We observe that  the finite groups $N_L(T)/T$ and  $N_{E_6}(T)/T$ have isomorphic
2-Sylow subgroups (of order $2^7$). By
\cite[Lemme 3.a]{G1}, it follows  that
there exists a finite odd degree extension
$K'/K$ such that  $\alpha_{K'}$ belongs to the image
of $H^1(K',L) \to H^1(K',E_6)$. Hence the class
$\alpha_{K'}$ is isotropic with respect to the root $\alpha_6$.
By the list of Witt--Tits indices in~\cite{tits},
the class $\alpha_{K'}$ is isotropic with respect to  $I$.
So $\alpha_{K'}$ belongs to the image
of $H^1(K',L_I) \cong H^1(K',P_I) \to H^1(K',E_6)$.

It thus remains to show that
\[ \ed(L_I, 2) \leq 3 \, . \]
To prove this inequality, we need an explicit description
of the group $L_I$.  Recall that there is a natural inclusion
$\mu_2 \times \mu_2 =C(DL_I) \subset T_I=\bbG_m \times \bbG_m$;
see \cite[Proof of Proposition 14.a]{CP}.
Hence we have the  following commutative exact diagram
$$
\begin{CD}
&& 1 && 1 \\
&& @VVV @VVV && \\
0 @>>> C(DL_I) =\mu_2 \times \mu_2 @>>> DL_I=\Spin_8 @>>> {\rm PSO}_8 @>>> 1 \\
 && @VVV @VVV  \mid \mid \\
0 @>>> T_I= \bbG_m \times \bbG_m  @>>> L_I @>>> {\rm PSO}_8  @>>> 1 .\\
 && @V{(2,2)}VV @VVV  \\
&& \bbG_m \times \bbG_m &=& \bbG_m \times \bbG_m  \\
 && @VVV @VVV  \\
 && 1 && 1\\
\end{CD}
$$
Taking Galois cohomology of the right square over
a field $F/k$, we obtain the following
commutative exact diagram of pointed sets
$$
\begin{CD}
&& H^1(F, \Spin_8) @>>> H^1(F, {\rm PSO}_8) \\
 && @VVV \mid \mid &&  \\
  1  @>>> H^1(F, L_I) @>>> H^1(F, {\rm PSO}_8)  .\\
  && @VVV  \\
&& 1    \\
 \end{CD}
$$
So $H^1(F, L_I)= {\rm Im}\Bigl( H^1(F, \Spin_8) \to
H^1(F, {\rm PSO}_8) \Bigr)$. It is well known that $H^1(F, L_I)$
classifies the similarity classes of $8$-dimensional quadratic $F$-forms
whose class belong to $I^3(F)$; cf., e.g.,~\cite[pp. 409 and 437]{boi}.
By the Arason-Pfister Theorem,
every  $8$-dimensional quadratic form $q \in I^3(F)$ is similar
to a 3-fold Pfister form $\langle \langle a,b,c \rangle \rangle$.
Thus the similarity class of $q$ is defined over $k(a, b, c)$. This
shows that $\ed(L_I,2) \leq 3$, as claimed.
\end{proof}

\begin{remark} \label{rem.E6} One can show that for 
every $\alpha \in H^1(K, E_6)$ there is an odd degree field
extension $L/K$ such that $\alpha_L$ lies in the image of 
the natural map $H^1(L, G_2) \to H^1(L, E_6)$; 
see~\cite[Exercise 22.9]{GMS}.
Since $\ed(G_2) = 3$, this leads to an alternative proof 
of the inequality $\ed(E_6; 2) \le 3$.
\end{remark}

\begin{example} \label{ex.E_7}
If $\Char(k) \ne 3$ then $\ed(E_7; 3) = 3$.
\end{example}

Here $E_7$ denotes the simply connected simple groups of type $E_7$.
By abuse of notation we will sometimes also
write $E_7$ for the Dynkin diagram of $E_7$.

\begin{proof}
By \cite[Table III]{griess}, $E_7$ has a unique non-toral
subgroup $H$ isomorphic to $(\bbZ/ 3 \bbZ)^5$.
To compute the rank of its centralizer, we make use
of its Witt--Tits index $I(H) \subset E_7$ \cite[section 3]{GP}, where
$E_7$ is the following Dynkin diagram
$$
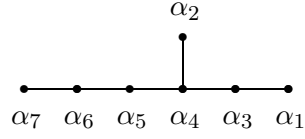
\begin{array}{ll}
\begin{picture}(105,45)
\put(00,02){\line(1,0){100}}
\put(60,02){\line(0,1){20}}
\put(00,02){\circle*{3}}
\put(20,02){\circle*{3}}
\put(40,02){\circle*{3}}
\put(60,02){\circle*{3}}
\put(80,02){\circle*{3}}
\put(100,02){\circle*{3}}
\put(60,22){\circle*{3}}
\put(20,02){\circle*{3}}
\put(60,02){\circle*{3}}
\put(80,02){\circle*{3}}
\put(100,02){\circle*{3}}
\put(-5,-11){$\alpha_7$}
\put(15,-11){$\alpha_6$}
\put(35,-11){$\alpha_5$}
\put(55,-11){$\alpha_4$}
\put(75,-11){$\alpha_3$}
\put(95,-11){$\alpha_1$}
\put(55,30){$\alpha_2$}
\end{picture}
\end{array}
$$

\medskip
\bigskip

\noindent Set $I= E_7 \setminus \{ \alpha_7\}$ and let
$P_I$ be the standard parabolic subgroup.
Denote by $L_I = Z_G(T_I)$ its standard Levi subgroup.
Then $DL_I=  E_6$, where $E_6$ denotes a simply connected group
of type $E_6$, Since  $E_6$ has
a  maximal non-toral  3-elementary abelian subgroup
of rank 4 (see \cite[Table III]{griess}),
we may assume that $H \subset E_6$.
Moreover, $C_{E_6}(H)$ is finite, so $H$ is irreducible in $L_I$.
It follows that $P_I$ is a minimal parabolic subgroup of $E_7$
containing $H$ and thus the Witt-Tits index of $H$ is $I$.
By \cite[Proposition 3.19]{GP}, the group $C_{E_7}(H)$ is of rank $1$.
Theorem~\ref{thm.main}(b) now tells us that $\ed(E_7; 3) \ge 4-1= 3$.

To prove the opposite inequality, consider $\alpha \in H^1(K,E_7)$.
By \cite[Example 3.5]{Ga}, the natural map
$$
H^1(K, E_6 \rtimes \mu_4) \to H^1(K,E_7)
$$
is surjective. Here, once again,
$E_6$ stands for the simply connected group of type $E_6$.
It follows that there exists a quartic extension $K'/K$
such that $\alpha_{K'}$ admits reduction of structure to $E_6$
(i.e., lies in the image of the map
 $H^1(K', E_6 ) \to H^1(K',E_7)$). Thus we may assume without
loss of generality that $\alpha$ comes from $E_6$.
Now recall that the natural map $H^1(K, F_4 \rtimes \mu_3) \to H^1(K,E_6)$
is surjective (see \cite[Example 3.5]{Ga}); here $\mu_3=C(E_6)$.
Thus there exists a $\beta \in H^1(K, F_4 \rtimes \mu_3)$ mapping to
$\alpha$ by the composite map
$$
H^1(K, F_4 \rtimes \mu_3) \to H^1(K,E_6) \to H^1(K,E_7).
$$
We claim that $\alpha \in H^1(K, E_7)$
admits further reduction of structure to $F_4$
(i.e., $\alpha$ depends only of the $F_4$-component of $\beta$).
If this claim is established the desired inequality
$\ed(\alpha;3) \leq 3$ will immediately follow from Example~\ref{ex.F_4}.

To prove the claim we again view $E_6$ inside $E_7$ as $E_6=DL_I$. We have
$L_I= C_{E_7}(T_I)=DL_I.T_I$, where $T_I=\bbG_m$ is the standard
torus associated to $I$ and $C(L_I)= T_I . C(E_7)$. Since $C(E_7)=\mu_2$,
it follows that $\mu_3=C(E_6) \subset T_I \subset L_I$.
We consider the commutative diagram of pairings
$$
\begin{CD}
H^1(K, \mu_3)& \times&  H^1(K,E_6) @>>> H^1(K,E_6) \\
@VVV @VVV @VVV \\
H^1(K, T_I) &\times & H^1(K, L_I) @>>> H^1(K, L_I) .\\
\end{CD}
$$
 From the vanishing of $H^1(K,T_I)$, it follows that the map
$H^1(K,E_6) \to H^1(K, L_I)$ is $H^1(K, \mu_3)$-invariant.
A fortiori, the image of the map
$$
H^1(K, F_4)  \to H^1(K,E_6) \to H^1(K, L_I)
$$
is $H^1(K, \mu_3)$-invariant. We conclude that the image of
$\beta \in H^1(F_4 \times \mu_3)$ in $H^1(K,L_I)$
depends only of its $F_4$ component, as claimed.
\end{proof}

\begin{remark} \label{rem.ex3-4}
One can show that if $G \to G'$ is a central isogeny of degree $d$
then $\ed(G; p) = \ed(G'; p)$ for any prime $p$ not dividing $d$.
In particular, the equalities $\ed(E_6; 2) = 3$ and $\ed(E_7; 3) = 3$
are valid for adjoint $E_6$ and $E_7$ as well. We leave the details
of this argument to the reader.
\end{remark}

\section*{Acknowledgements}
We would like to thank J.-P. Serre for sharing his conjecture
with us and for numerous insights related to the subject
matter of this paper. We are also grateful to M. Brion and O. Gabber
for helpful comments.

The second author would like to thank Laboratoire de math\'ematiques
at the Universit\'e Paris-Sud for their hospitality during
the week of June 19--25, 2006, when the work on this collaborative
project began.

\end{document}